\documentclass[11pt]{article}

\usepackage{amsmath}
\usepackage{amsthm}
\usepackage{amssymb}
\makeatletter

\newtheorem{thm}{Theorem}[section]
\newtheorem{cor}[thm]{Corollary}
\newtheorem{prob}[thm]{Problem}
\date{}

\let\originalleft\left
\let\originalright\right
\renewcommand{\left}{\mathopen{}\mathclose\bgroup\originalleft}
\renewcommand{\right}{\aftergroup\egroup\originalright}
\makeatletter

\let\Url@force@Tilde\UrlTildeSpecial
\makeatother

\makeatother

\begin{document}

%SK \title{Extremal problem for integer sparse recovery}
\title{An extremal problem for integer sparse recovery}

\author{Sergei Konyagin\thanks{Laboratory ``Multivariate approximation
and applications'', Department of Mechanics and Mathematics,
Moscow State University, Moscow, Russia; Steklov
Mathematical Institute of Russian Academy of Sciences, Moscow, Russia.
Email: konyagin@mi-ras.ru. Research supported
by the grant of the Government of the
Russian Federation (project 14.W03.31.0031).}
\and
Benny Sudakov\thanks{Department of Mathematics, ETH Zurich, Switzerland. Email: benjamin.sudakov@math.ethz.ch.
Research supported in part by SNSF grant 200021-175573.}
}

\maketitle
\global\long\def\cA{\mathcal{A}}
\global\long\def\Ex{\mathbb{E}}
\global\long\def\RR{\mathbb{R}}
\global\long\def\E{\mathbb{E}}
\global\long\def\cF{\mathcal{F}}
\global\long\def\Var{\operatorname{Var}}
\global\long\def\NN{\mathbb{N}}
\global\long\def\ZZ{\mathbb{Z}}
\global\long\def\eps{\varepsilon}
\global\long\def\one{\boldsymbol{1}}
\global\long\def\range#1{\left[#1\right]}
\global\long\def\Bin{\operatorname{Bin}}
\global\long\def\cT{\mathcal{T}}
\global\long\def\floor#1{\left\lfloor #1\right\rfloor }
\global\long\def\ceil#1{\left\lceil #1\right\rceil }

\begin{abstract}
Motivated by  the problem of integer sparse recovery we study the following question.
%SK2 Let $A$ be an $m \times d$ integer matrix whose entries have absolute value at most $k$.
Let $A$ be an $m \times d$ integer matrix whose entries are in absolute value at most $k$.
How large can be $d=d(m,k)$ if all $m \times m$ submatrices of $A$ are non-degenerate?
We obtain new upper and lower bounds on $d$ and answer  
a special case of the problem by Brass, Moser and Pach on covering $m$-dimensional $k \times \cdots\times k$ grid by
linear subspaces. 
\end{abstract}

\section{Introduction}

Compressed sensing is a relatively new mathematical paradigm that shows a small number of linear
measurements are enough to efficiently reconstruct a large dimensional signal under the assumption
that the signal is \textit{sparse} (see, e.g., \cite{FR} and its references). That is, given a signal ${\bf{x}}\in {\mathbb R}^d$, the
%SK2 goal is to accurately reconstruct $\bf{x}$ from its noisy measurements $\bf{b} = A\bf{x} + \bf{e}$.
goal is to accurately reconstruct $\bf{x}$ from its noisy measurements ${\bf b} = A{\bf x} + {\bf e}$.
Here, $A$ is an underdetermined matrix $A\in \mathbb{R}^{m\times d}$, where $m$ is much smaller
than $d$,
and ${\bf{e}} \in {\mathbb R}^m$ is a vector modeling noise in the system. Since the system is highly 
underdetermined, it is ill-posed until one imposes additional constraints, such as the signal 
%SK2 $\bf{x}$ obeying a sparsity constrain. We say $\bf{x}$ is {\it $s$-sparse} when it has at most 
$\bf{x}$ obeying a sparsity constraint. We say $\bf{x}$ is {\it $s$-sparse} when it has at most 
%SK2 $s$ nonzero entries. Clearly, any matrix $A$ that is one-to-one on $2s$-sparse signals will allow 
%SK2 Maybe, it is reasonable to revise this sentence 
$s$ nonzero entries. Clearly, any matrix $A$ that is one-to-one on $s$-sparse signals will allow 
reconstruction in the noiseless case when $\bf{e}=0$. However, compressed sensing seeks the ability to reconstruct efficiently and robustly even when one allows presence of noise.
Motivated by this problem Fukshansky, Nadel and Sudakov~\cite{FKS} considered the following
extremal problem, which is of independent interest.

\begin{prob}
\label{P1}
%SK2 Given an integers $k,m$ what is the maximum integer $d$ such that 
Given integers $k,m$, what is the maximum integer $d$ such that 
there exists $m \times d$ matrix $A$ with integer entries satisfying $|a_{ij}|\leq k$
%SK2 such that all $m \times m$ submatrices of $A$ are non-degenerate.
such that all $m \times m$ submatrices of $A$ are non-degenerate?
\end{prob}
To see the connection of this question with integer sparse recovery let $s \leq m/2$ and consider
$s$-sparse signal ${\bf{x}} \in \mathbb{Z}^d$. We denote by $\|\bf{b}\|$ the Euclidean norm
of a vector ${\bf b}=(b_1,\dots,b_m)\in\mathbb{R}^m$ and by $\|\bf{b}\|_\infty$ its $l_\infty$-norm:
$\|{\bf b}\|_\infty=\max_{i=1,\dots,m}|b_i|$.
Suppose we wish to decode ${\bf{x}}$ from the noisy measurements ${\bf b} = A{\bf x} + {\bf e}$
where $\|{\bf{e}}\|_\infty < \frac{1}{2}$ (in particular, this
holds if $\|{\bf{e}}\| < \frac{1}{2}$). Note that by definition of
%SK2 matrix $A$ we have that for any $m$-sparse vector $\bf{z}$,
matrix $A$ we have for any $m$-sparse integer non-zero vector $\bf{z}$ that
%SK2 $A{\bf z} \not =0$ and therefore being integer vector has
$A{\bf z} \not ={\bf 0}$ and therefore, being an integer vector, $A{\bf z}$ has
$l_\infty$-norm at least one. So to decode $\bf{x}$
%SK2 we can select the $s$-sparse signal ${\bf y}\in \mathbb{Z}^d$ minimizing $\|{\bf b} - A{\bf y}\|$.
we can select the $s$-sparse signal ${\bf y}\in \mathbb{Z}^d$ minimizing $\|{\bf b} - A{\bf y}\|_\infty$.
%SK2 Then since $\bf{x}$ satisfies  ${\|{\bf b}- A{\bf x}\|_ \infty= \|{\bf e}\|_\infty} < \frac{1}{2}$,
Then, since $\bf{x}$ satisfies  ${\|{\bf b}- A{\bf x}\|_ \infty= \|{\bf e}\|_\infty} < \frac{1}{2}$,
it must be that the decoded vector $\bf{y}$ satisfies
this inequality as well. Therefore,
$\|A{\bf y} - A{\bf x}\|_\infty \leq \|{\bf b} - A{\bf y}\|_\infty + \|{\bf b} - A{\bf x} \|_\infty< 1$ 
is an $m$-sparse vector, which guarantees that
%SK3 ${\bf y}={\bf x}$ so our decoding was successful. Note that if instead of error $1/2$ we want to
%SK3 allow larger error $C$ we can simply multiply all entries of $A$ by a factor of $2C$.
${\bf y}={\bf x}$ so our decoding was successful. Note that if instead of error $<1/2$ we want to
allow error $<C$ we can simply multiply all entries of $A$ by a factor of $2C$.

Fukshansky, Nadel and Sudakov~\cite{FKS} showed that a matrix $A\in {\mathbb R}^{m\times d}$ with integer
entries $|a_{ij}|\leq k$ and all $m \times m$ submatrices having full rank must satisfy
$d=O(k^2m)$. They also proved that such matrices exists when
$d =\Omega(\sqrt{k} m)$. Their upper bound was improved by Konyagin~\cite{K} who showed that
$d$ must have order at most  $O(k (\log k) \,m)$ (all logarithms here and later in the paper
are in base $e$) for $m\geq\log k$ and at most $O(k^{m/(m-1)}\,m^2)$ for $2\leq m<\log k$.
Improving these results further, in this paper we obtain the following new upper bound.

\begin{thm}
\label{th1}
%SK2 Let $A$ be an $m \times d$ integer matrix such that $|a_{ij}|\leq k$, $m\ge\log k$,
Let $A$ be an $m \times d$ integer matrix such that $|a_{ij}|\leq k$,
and all $m \times m$
%SK2 submatrices of $A$ have a full rank. If $k$ is sufficiently large, then
submatrices of $A$ have full rank. If $k$ is sufficiently large, then
$d \leq 100k \sqrt{\log k}\, m$ for $m\geq\log k$ and $d\leq 400 k^{m/(m-1)}\,m^{3/2}$
for $2\leq m < \log k$.
\end{thm}

The lower bound construction for Problem \ref{P1} uses random matrices and is based on a deep result of Bourgain, Vu and Wood~\cite{BVW}  which estimates the probability that a random $m \times m$ matrix with integer entries from $[-k, k]$ is singular.
It is expected that their result is not tight and the probability of singularity for such matrix
%SK has order $k^{(1-o(1))m}$. If this is the case
has order $k^{-(1-o(1))m}$ as $k\to\infty$. If this is the case
then $m \times d$ matrices, satisfying  Problem \ref{P1}, exist for $d$ close to $km$.
This suggests that our new bound for $m\geq\log k$ is not far from being optimal.

On the other hand, we get the following result.

\begin{thm}
\label{th2}
Let $k\in\NN$, $m\in\NN$, $m\ge2$, and
$$m < d \leq\max(k+1, k^{m/(m-1)}/2).$$
Then there is an $m \times d$ integer matrix $A$ such that $|a_{ij}|\leq k$
%SK2 and all $m \times m$ submatrices of $A$ have a full rank.
and all $m \times m$ submatrices of $A$ have full rank.
\end{thm}

We observe that this theorem improves the lower bound from~\cite{FKS}
for $m=o(\sqrt k)$. Moreover, the existence of required matrices in~\cite{FKS}
was proven by using probabilistic arguments, but the matrices in Theorem~\ref{th2}
are explicit and easily computable. Also we notice that upper and lower estimate for
$d$ in Theorems \ref{th1} and \ref{th2} differ by a factor $O(m^{3/2})$
depending on $m$ only.

%SK2 Our result can be also used to answer a special case of the problem 
Our result can be also used to solve a special case of a problem 
by Brass, Moser and Pach. In \cite{BMP} (Chapter 10.2, Problem 6) they 
%SK2 asked, what is the minimum number $M$ of $s$-dimensional linear subspaces  
asked: what is the minimum number $M$ of $s$-dimensional linear subspaces  
necessary to cover $m$-dimensional $k \times \cdots\times k$ grid 
%SK2 $K=\{x\in\ZZ^m:\|x\|_\infty\le k\}$ Balko, Cibulko and Valtr~\cite{BCV} 
$K=\{{\bf x}\in\ZZ^m:\|{\bf x}\|_\infty\le k\}$? Balko, Cibulko and Valtr~\cite{BCV} 
studied this problem and obtain upper and lower bounds for $M$.
In particular in the case when $s=m-1$ they proved that
$$k^{m/(m-1)-o(1)} \le M\le C_m k^{m/(m-1)}.$$
Using Theorem~\ref{th2} we obtain a new lower bound which is tight up to a constant factor.

\begin{cor}
\label{cor1}
For $k\geq m\geq 2$ we have
$$M \geq  k^{m/(m-1)}/(2m-2).$$
\end{cor}

Indeed, suppose that we cover $K$ by $M$ hypersubspaces $P_1,\dots,P_M$.
We consider the columns of the matrix $A$ constructed in Theorem~\ref{th2}.
Since any $m$ of them are linearly independent, every subspace $P_i$
contains at most $m-1$ of these columns. Thus, $d\le M(m-1)$, and the corollary
follows.

%SK2
%SK3 For matrices $A$ constructed in Theorem \ref{th2}, Ryutin \cite{R} suggested
%SK3 an effective algorithm for recovery of a sparse vector $\bf x$ via $A\bf x$.
For matrices $A$ constructed in Theorem \ref{th2}, Ryutin \cite{R} suggested
an effective algorithm for recovery of a sparse vector $\bf x$ via $A\bf x$.
See also \cite{R2}.

%\vspace{0.2cm}
%\noindent

\section{Proof of Theorem \ref{th1}}

Let $t=[\log k]$ for $m\geq\log k$ and
$t=m$ for $2\leq m < \log k$. In the first case we suppose that
$d > 100k \sqrt{t} m$ and in the second case we suppose that
%SK2 $d > 400 k^{m/(m-1)}\,m^{3/2}$.
$d > 400 k^{m/(m-1)}\,\sqrt{t}m$.
Let $v_1, \dots, v_t$ be the first $t$ rows of the matrix $A$.
Take $\Lambda=9$ if $m\geq\log k$ and $\Lambda = [25 k^{1/(m-1)}]$ otherwise. Given a
%SKmathbold vector of integer coefficients $\mathbold{\lambda}=(\lambda_1, \ldots, \lambda_t)$
vector of integer coefficients ${\bf{\lambda}}=(\lambda_1, \ldots, \lambda_t)$
%SKmathbold such that $0 \leq \lambda_i \leq 9$ denote by $v_{\mathbold{\lambda}}$
such that $0 \leq \lambda_i \leq \Lambda$ denote by $v_{\bf\lambda}$
a linear combination $\sum_i \lambda_i v_i$.
%SKmathbold Our goal is to find two combinations $\mathbold{\lambda} \not = \mathbold{\lambda}'$
Our goal is to find two combinations ${\bf\lambda} \not = {\bf\lambda}'$
%SKmathbold such that corresponding vectors $v_{\mathbold{\lambda}}$ and $v_{\mathbold{\lambda}'}$
such that corresponding vectors $v_{\bf\lambda}$ and $v_{\bf\lambda}'$
agree on at least $m$ coordinates.
This will show that a linear combination of first $t$ rows of matrix $A$ with coefficients
%SKmathbold $\mathbold{\lambda}-\mathbold{\lambda}'\not =\bf{0}$ has at least $m$ zeros
${\bf\lambda}-{\bf\lambda}'\not =\bf{0}$ has at least $m$ zeros
and therefore the $m \times m$ submatrix of $A$
whose columns correspond to these zeros is degenerate, since its first $t$ rows are
linearly dependent.

Consider ${\bf\lambda}$ chosen uniformly at random out of $(\Lambda+1)^t$
%SK2 possible vectors and look on a value of a fixed coordinate $j$ of the vector
possible vectors and look at a value of a fixed coordinate $j$ of the vector
$v_{\bf\lambda}$.
This value is a random variable $X$ which is a sum of the $t$
%SK2 independent random variable $X_i$, where $X_i$ is a value of the $j$-th coordinate of
independent random variables $X_i$, where $X_i$ is a value of the $j$-th coordinate of
$\lambda_i v_i$. Since $|a_{ij}|\leq k$, we have that $|X_i| \leq \Lambda k$ and therefore its
variance $\Var(X_i) \leq \mathbb{E}(X_i^2) \leq \Lambda^2k^2$. This implies that
$\Var(X)=\sum_i\Var(X_i) \leq \Lambda^2k^2t$. Thus, by Chebyshev's inequality, with probability
%SK2 at least $3/4$ the value of $X$ belongs to an interval $I$ 
at least $3/4$, the value of $X$ belongs to an interval $I$ 
of length $4\sqrt{\Var(X)} \leq 4\Lambda k\sqrt{t}$.
%SK2 Hence there are at least $0.75 \cdot (\Lambda+1)^t$ linear combinations $v_{\bf\lambda}$
Hence there are at least $0.75 \cdot (\Lambda+1)^t$ values of $\bf\lambda$ giving
a linear combination $v_{\bf\lambda}$
whose $j$-th coordinate belongs to $I$. For every integer $s$ let $h_j(s)$ be the number of
%SK2 $v_{\bf\lambda}$ whose $j$-th coordinate is $s$ and let $h_j$ be the number
values for $\bf\lambda$ giving linear combinations $v_{\bf\lambda}$ whose $j$-th coordinate 
is $s$, and let $h_j$ be the number
%SKmathbold of ordered pairs $\mathbold{\lambda} \not = \mathbold{\lambda}'$
of ordered pairs ${\bf\lambda} \not = {\bf\lambda}'$
%SKmathbold such that $v_{\mathbold{\lambda}}$ and $v_{\mathbold{\lambda}'}$ agree on $j$-th coordinate.
such that $v_{\bf\lambda}$ and $v_{\bf\lambda}'$ agree on $j$-th coordinate.
By definition $0.75 \cdot (\Lambda+1)^t \leq \sum_{s \in I}  h_j(s) \leq (\Lambda+1)^t$
and $h_j=\sum_s h_j(s)(h_j(s)-1)$.

%SK2 If $m\geq\log k$ and $\Lambda=9$ then using a Cauchy-Schwarz inequality,
If $m\geq\log k$, in which case $\Lambda=9$ and $t=[\log k]$, then using the Cauchy-Schwarz inequality,
together with the facts that
%SK2 $10^t>k^2\gg k\sqrt{\log k}=k\sqrt{t}$ for sufficiently large $k$,  we have
$10^t>k^2\geq k\sqrt{\log k}=k\sqrt{t}$ for sufficiently large $k$ and the interval
$I$ contains at most $4\Lambda k\sqrt{t}+1$ integer points, we have
\begin{eqnarray*}
h_j&=&\sum_s h_j(s)(h_j(s)-1) = \sum_s  h^2_j(s)-\sum_s  h_j(s) \geq \sum_{s \in I}  h^2_j(s)-10^t\\
%SK2 &\geq& \frac{1}{4 \Lambda k\sqrt{t}}\left(\sum_{s \in I}  h_j(s)\right)^2-10^t
&\geq& \frac{1}{4 \Lambda k\sqrt{t}+1}\left(\sum_{s \in I}  h_j(s)\right)^2-10^t
\geq \frac{1}{40k\sqrt{t}} \big(0.75 \cdot 10^t\big)^2-10^t \\
&\geq&  \frac{1}{80k\sqrt{t}}10^{2t}-10^t \geq \frac{1}{100k\sqrt{t}}10^{2t}\,.
\end{eqnarray*}
%SKmathbold Since the number of ordered pairs $\mathbold{\lambda}\not = \mathbold{\lambda}'$
Since the number of ordered pairs ${\bf\lambda}\not = {\bf\lambda}'$
is at most $10^{2t}$ and the number of coordinates $j$ is $d$, by averaging we obtain that
there is a pair ${\bf\lambda}\not = {\bf\lambda}'$
which agrees on at least
%SK2 $$\frac{\sum_{j=1}^d h_j}{10^{2t}} \geq \frac{d}{100k\sqrt{t}} \geq m$$
$$\frac{\sum_{j=1}^d h_j}{10^{2t}} \geq \frac{d}{100k\sqrt{t}}$$
%SK2 coordinates. As we explain above, this implies that $A$ has an $m \times m$
%SK2 degenerate submatrix.
coordinates. If the pair ${\bf\lambda}\not = {\bf\lambda}'$ agrees on 
$\frac{d}{100k\sqrt{t}} \geq m$ coordinates, then as 
we explain above, this implies that $A$ has an $m \times m$
degenerate submatrix. Thus, we must have 
$$d\leq 100 km\sqrt t= 100 km\sqrt{\log k}$$
when $m\geq\log k$.

Now we consider the case $2\leq m < \log k$. Then, due to the inequality
$$\frac{(\Lambda+1)^m}{8\Lambda k\sqrt{m}} \geq
\frac{25^m k^{m/(m-1)}}{200k^{m/(m-1)}\sqrt{m}} > 2,$$
we have

\begin{eqnarray*}
h_j&=&\sum_s h_j(s)(h_j(s)-1) = \sum_s  h^2_j(s)-\sum_s  h_j(s) \geq \sum_{s \in I}  h^2_j(s)-
(\Lambda+1)^m\\
&\geq& \frac{1}{4\Lambda k\sqrt{m}+1}\left(\sum_{s \in I}  h_j(s)\right)^2-(\Lambda+1)^m\\
&\geq& \frac{1}{4\Lambda k\sqrt{m}+1} \big(0.75 \cdot (\Lambda+1)^m\big)^2-(\Lambda+1)^m \\
&\geq& (\Lambda+1)^m\left(\frac{(\Lambda+1)^m}{8\Lambda k\sqrt{m}} - 1\right)
\geq (\Lambda+1)^m \frac{(\Lambda+1)^m}{16\Lambda k\sqrt{m}}.
\end{eqnarray*}
Since the number of ordered pairs ${\bf\lambda}\not = {\bf\lambda}'$
is at most $(\Lambda+1)^{2m}$ and the number of coordinates $j$ is $d$, by averaging
we obtain that there is a pair ${\bf\lambda}\not = {\bf\lambda}'$
which agrees on at least
%SK2 $$\frac{\sum_{j=1}^d h_j}{(\Lambda+1)^{2m}} \geq \frac{d}{16\Lambda k\sqrt{m}} \geq m$$
$$\frac{\sum_{j=1}^d h_j}{(\Lambda+1)^{2m}} \geq \frac{d}{16\Lambda k\sqrt{m}}$$
%SK2 coordinates as required. This completes the proof of the theorem. \hfill $\Box$
coordinates. If the pair ${\bf\lambda}\not = {\bf\lambda}'$ agrees on
$\frac{d}{16\Lambda k\sqrt{m}} \geq m$ coordinates, then as 
we explain above, this implies that $A$ has an $m \times m$
degenerate submatrix. Thus, we must have 
$$d\leq 16\Lambda k m\sqrt m\leq 400 k^{m/(m-1)}\,m^{3/2}$$
when $2\leq m<\log k$. This completes the proof of the theorem. \hfill $\Box$

\section{Proof of Theorem \ref{th2}}

%SK2 We have to construct required $m\times d$ matrices
We have to construct the required $m\times d$ matrices
%SK2 with $d\geq k+1$ provided that $k\geq m$ and
with $d\geq k+1$, provided that $k\geq m$, and
%SK2 $d\geq k^{m/(m-1)}/2$ provided that $k^{m/(m-1)}/2>m$.
with $d\geq k^{m/(m-1)}/2$, provided that $k^{m/(m-1)}/2>m$.

First we will construct an $m\times d$ matrix with $d\geq k+1$.
%SK2 There exists an odd prime $d$ with $k+1\leq  d\leq 2k+1$.
By Bertrand's postulate, there exists an odd prime $d$ with $k+1\leq  d\leq 2k+1$.
%SK2 We define the matrix $A$ by taking $a_{i,j}\equiv j^{i-1} (\mod d)$
We define the matrix $A$ by taking $a_{i,j}\equiv j^{i-1} (\bmod d)$
with $|a_{i,j}|\le(d-1)/2\le k$. Considering the matrix $A$ modulo $d$ we find
that any submatrix of $A$ of size $m\times m$ is a Vandermonde matrix modulo
$d$. Hence, its determinant is not zero modulo $d$. This implies that this
%SK2 submatrix has a full rank.
submatrix has full rank.

Next we will construct an $m\times d$ matrix with $d\geq k^{m/(m-1)}/2$.
%SK2 We can consider that $k^{m/(m-1)}/2>k+1$ and, in particular, $k\geq 3$.
We can assume that $k^{m/(m-1)}/2>k+1$ and, in particular, $k\geq 3$.
%SK2 There exists a prime $d$ with $k^{m/(m-1)}/2\leq  d < k^{m/(m-1)}$.
By Bertrand's postulate, there exists a prime $d$ with $k^{m/(m-1)}/2\leq  d < k^{m/(m-1)}$.
For $u\in\RR$ we denote by $\|u\|$ the distance from $u$ to the nearest
integer. Note that if we can multiply $u$ by some integer $\ell$ such that $\ell u$ is also an integer
then the distance of $\ell u$ to the closest multiple of $\ell$ is exactly $\ell \|u\|$. Therefore we can find and 
%SK3 integer $z$ such that $z=\ell u \mod \ell$ and $|z|=\ell \|u\|$
integer $z$ such that $z\equiv\ell u(\bmod \ell)$ and $|z|=\ell \|u\|$.

We first consider the $m\times d$ matrix $A'$ with entries
$a_{i,j}'= j^{i-1}$. Again, the determinant of any $m\times m$
submatrix of $A'$ is not zero modulo $d$. The idea is to multiply the
columns of $A'$ by appropriate integers not divisible by $d$ such that when we
replace all entries by their residues modulo $d$, the absolute
values of this residues will  be bounded by $k$. Clearly, the operation of multiplying 
%SK3 by integers not divisible by $d$ and taking residues$\mod d $ preserve the
by integers not divisible by a prime $d$ and taking residues$\mod d $ preserve the
%SK3 property of submatrices having full rank.
property that determinants of all $m\times m$ submatrices are not divisble by $d$.

Using Dirichlet's theorem on simultaneous approximations (see, e.g., \cite{S}, Chapter 2, Theorem 1A), we find that for
every $j=1,\dots,d$ there is a positive integer $l_j<d$
such that $\|l_j j^{i-1}/d\|\leq d^{-1/m}$ for $i=1,\dots,m$.
Hence, for any $i$ there is an integer $a_{i,j}$ such that
%SK2 $a_{i,j} \equiv l_j j^{i-1}(\mod d)$ and
$a_{i,j} \equiv l_j j^{i-1}(\bmod d)$ and
$|a_{i,j}|\le d^{1-1/m}\leq k$ as required. This completes
the proof of Theorem~\ref{th2}.\hfill $\Box$

\vspace{0.3cm}

\noindent
{\bf Acknowledgement.}\, We thank the anonymous referees for many useful suggestions.

\end{document}